\documentclass[12pt,centertags,oneside]{amsart}

\usepackage{amsmath,amstext,amsthm,amscd,typearea,hyperref,color}
\usepackage{amssymb}
\usepackage{a4wide}
\usepackage[mathscr]{eucal}
\usepackage{mathrsfs}
\usepackage{typearea}
\usepackage{charter}
\usepackage{pdfsync}
\usepackage[a4paper,width=16.2cm,top=3cm,bottom=3cm]{geometry}

\numberwithin{equation}{section}

\newtheorem{theorem}{Theorem}[section]

\newtheorem{problem}[theorem]{Problem}

\newcommand{\cali}[1]{\mathscr{#1}}

\newcommand{\ddc}{{dd^c}}

\newcommand{\ddbar}{{\partial\overline\partial}}

\newcommand{\sing}{{\rm sing}}
\newcommand{\Sing}{{\rm Sing}}

\newcommand{\Ac}{\cali{A}}
\newcommand{\Bc}{\cali{B}}

\newcommand{\Dc}{\cali{D}}

\newcommand{\Fc}{\cali{F}}

\newcommand{\Hc}{\cali{H}}

\renewcommand{\Mc}{\cali{M}}

\newcommand{\Uc}{\cali{U}}

\newcommand{\C}{\mathbb{C}}
\newcommand{\D}{\mathbb{D}}

\newcommand{\T}{\mathbb{T}}

\renewcommand{\P}{\mathbb{P}}

\title{Some open problems on holomorphic foliation theory}

\author{Tien-Cuong Dinh}
\address{Department of Mathematics, National University 
of Singapore, 10 Lower Kent Ridge Road, Singapore 119076.}
\email{matdtc@nus.edu.sg}

\author{Nessim Sibony}
\address{Laboratoire de Math\'ematiques d'Orsay, Univ. Paris-Sud, CNRS, Universit\'e
Paris-Saclay, 91405 Orsay, France
\endgraf and\endgraf 
Korea Institute for Advanced Study, 
Seoul, 130-722
South KOREA}

\email{Nessim.Sibony@math.u-psud.fr}

\date{July 10, 2018}

\begin{document}

\maketitle

\begin{abstract}
We present a list of open questions in the theory of holomorphic foliations, possibly with singularities. Some problems have been around for a while, others are very accessible.
\end{abstract}

\medskip\medskip

\noindent
{\bf Classification AMS 2010:} 37F75, 37A, 32U40, 30F15, 57R30.

 \medskip

\noindent
{\bf Keywords:} Singular foliation, $\ddc$-closed current, Nevanlinna current, unique ergodicity.

\section{Introduction} \label{introduction}

The study of differential equations on manifolds is at the center of several mathematical questions.
Although the subject is of classical nature, many important questions are unsolved. 
Here, we discuss some problems in the context of meromorphic vector fields in the complex projective space $\P^k$. Since the vector fields are not holomorphic, it is useful to rather consider foliations.

In the field of holomorphic foliations, even in $\P^2,$ the following are examples of open basic questions : 
Do leaves always cluster at singular points ? What is the topology of generic leaves ? Are leaves dense for generic foliations ?

To give a structure to the space of foliations of a given degree, we  use the following  point of view.
Let $\pi: \C^{k+1} \setminus \{0\}\rightarrow{\P}^k$
denote the canonical projection. Consider the foliation defined
in $ \C^{k+1}$ by a vector field
$$Z(z) := \sum_{j=0}^{k} F_j(z) \frac{\partial}{\partial z_j}$$
which is not identically zero.
Here, the $F_j$'s are homogeneous holomorphic polynomials of a fixed degree $d$ without non-trivial common divisor. The homogeneity 
permits to descend the integral curves of $Z$ in $\C^{k+1}$ to a foliation by Riemann surfaces $\Fc$ in $\P^k$.
Moreover, all holomorphic foliations in $\P^k$ can be obtained in this way.

So we have constructed the space $\Bc_d$ of holomorphic foliations of degree $d\geq 2$ in $\P^k$, which can be parametrized by 
a dense Zariski open set in some projective space $\P^N$. 
The number $d$ is {\it the geometric degree} of the foliation and can be interpreted as 
 the number of tangency points of leaves of $\Fc$ with a generic projective hyperplane in $\P^k$. 
 Locally, the foliation $\Fc$  is tangent to a family of $k-1$ holomorphic $(1,0)$-forms $(\alpha_j)_{1\leq j\leq k-1}$ whose values at a generic point are linearly independent.
To study the global behavior of leaves, it is useful to introduce currents directed by the foliation which are a version of the harmonic measures of Lucy Garnett \cite{G1983} for non-singular foliations. 

A positive closed  (resp. $\ddc$-closed) current $T$ of bi-dimension $(1,1)$ is {\it tangent to} (or {\it directed by}) the foliation $\Fc$  if $T\wedge \alpha_j=0$ for $1\leq j\leq k-1$. Recall that $\ddc:={i\over\pi}\ddbar$. 
Such currents play an important role in the ergodic theory of foliations, see \cite{DNS1, DNS2, DNS3, DS1, FS1, FS3, Nguyen1, Nguyen2}. In particular, the survey \cite{FS3} provides some background.

The point of view of currents, rather than measures, provides more structure (e.g.  the 
notion of intersection and density for currents) which has been at the center of recent progress in the field.

\section{Ergodicity for foliations} \label{s:unique}

In \cite{DNS1}, Dinh-Nguyen-Sibony proved the following geometrical ergodic theorem. Let $\Fc$ be a foliation, possibly singular, by Riemann surfaces on a compact K\"ahler manifold $M$. The manifold $M$ is endowed with a fixed K\"ahler form $\omega$.  Let $T$ be an extremal positive $\ddc$-closed current directed by the foliation $\Fc$. Assume for simplicity that $T$ has full mass on the hyperbolic leaves, i.e. leaves covered by the unit disc.

If $\phi_a:\D\to L$ is the universal covering of a leaf $L$ with $\phi_a(0)=a$, we consider for $r<1$ the {\it Nevanlinna current}
$$\tau_r^a:= {1\over T(r)} (\phi_a)_*\Big(\log^+{r\over |\xi|}\Big) = {1\over T(r)} \int_0^r {dt\over t} (\phi_a)_*[\D(t)].$$
Here, $[\D(t)]$ is the current of integration on the disc $\D(t)$ of radius $t$ centered at  0, $\log^+:=\max(\log,0)$ and $T(r)$ is the {\it Nevanlinna characteristic} for $\phi_a$ which is given by the formula
$$T(r) := \int_0^r{dt\over t} \int_{\D(t)} (\phi_a)^*(\omega).$$

\begin{theorem}[\cite{DNS1}] \label{t:DNS}
In the above situation, assume that all singularities of $\Fc$ are linearizable. Then, for $T$-almost every $a$, we have $\tau_r^a\to T$ as $r$ tends to $1$, in the sense of currents. More precisely, for any smooth test form $\theta$ of bi-degree $(1,1)$ on $\P^k$, we have 
$$\langle \tau^a_r,\theta\rangle \to \langle T,\theta\rangle \quad\text{as} \quad r\to1.$$
\end{theorem}

\begin{problem}
(i) Is Theorem \ref{t:DNS} still true with non-linearizable singularities ?
The question here is to estimate the Poincar\'e metric on leaves near a non-linearizable singularity.

(ii) Is there a positive $\ddc$-closed current $T$, directed by $\Fc$, extremal, for which one can estimate the speed of the convergence in Theorem \ref{t:DNS} for smooth observables $\theta$ ?
\end{problem}

For the second question, a good candidate for $T$  is a current of positive entropy as defined in \cite{DNS2}.

Instead of considering the weighted current $\tau^a_r$, we can consider an un-weighted current ({\it Ahlfors current})
$$\sigma_r^a:={1\over A(r)} (\phi_a)_*[\D(r)]  \qquad \text{with} \qquad  A(r):=\int_{\D(r)} \phi_a^*(\omega).$$

\begin{problem}
Do we have $\sigma_r^a\to T$ as $r$ tends to $1$ for $T$-almost every $a$ ?
\end{problem}

A consequence of Theorem \ref{t:DNS} is that for almost every leaf, with respect to the current $T,$
the Nevanlinna characteristic $T(r)$ for the leaf parametrized by $\phi_a$, satisfies the inequality
$$ T(r)\geq c \log{1\over(1-r)}$$
for $r$ close to 1, where $c$ is a positive constant depending on $a$.

\begin{problem}
  Is the above estimate true for every leaf ? In particular, is the above estimate true for every leaf, when there is no diffuse $\ddc$-closed current ?
\end{problem}

Recall that a current is {\it diffuse} if it gives no mass to analytic varieties.
There are dynamically non-trivial cases, where the only  $\ddc$-closed currents directed by the foliation 
are supported by varieties. Indeed, assume $\Fc$ is a foliation in $\P^2$ with only hyperbolic singularities 
and the line at infinity is  the only invariant curve for $\Fc$. Then the only 
positive $\ddc$-closed current of mass 1 directed by $\Fc$  is the closed current defined by the line at infinity  \cite{DS1}.

Recall that a singularity $a$ of $\Fc$ is hyperbolic when the quotient of the eigenvalues at $a$ of the linearized vector field defining $\Fc$ is non-real.
So, in the considered case, the convergence in Theorem \ref{t:DNS} is true for every $a$.  Hence, by \cite{DS1},  $T(r)$ is unbounded but no explicit estimate for the growth is known. 
For a compact lamination, without singularities, it follows from \cite{FS1} that the above lower bound for $T(r)$ holds for all leaves.

\section{Unique ergodicity for foliations} \label{s:Unique}
The following unique ergodicity result was completed recently, see \cite {DNS3, DS1, FS1}.
It can be extended to compact K\"ahler surfaces, but we mention here only the $\P^2$ case for simplicity.

 \begin{theorem} \label{th_main_P2_bis}
 Let $\Fc$ be a holomorphic foliation by Riemann surfaces in $\P^2$ whose singularities are all hyperbolic. Then
 either (1) $\Fc$ has no invariant algebraic curve and then it has a unique directed positive $\ddc$-closed  $(1,1)$-current of mass $1$; moreover, this unique current is not closed, or (2) $\Fc$ admits a finite number of invariant algebraic curves and any directed  positive $\ddc$-closed $(1,1)$-current is a combination of currents of integration on these curves; in particular, all directed positive $\ddc$-closed $(1,1)$-currents are closed.
 \end{theorem}

\begin{problem}
Assume that $\Fc$ is a foliation of degree $d\geq 2$ with non-hyperbolic singularities. Describe the
$\ddc$-closed positive currents, directed by the foliation. In which case diffuse positive closed 
currents directed by the foliation do exist ?
\end{problem}
When the singularities are hyperbolic, i.e. when the quotient of the eigenvalues (of the linearized vector field) is non-real, the only local positive closed current near the singularity are
supported by the separatrices.  This is due to the contracting holonomy. 

When the quotient of the eigenvalues (of the linearized vector field) at a non degenerate singularity is real, then the space of local positive closed currents directed by the local foliation is quite big and easy to describe. However, the existence of global currents  is unclear. It is in particular unclear whether, one could have for a foliation in $\P^2$, both diffuse positive closed currents and diffuse positive $\ddc$-closed currents which are not closed (extremal $\ddc$-closed, non-closed, currents are necessarily diffuse). 

Let $\Bc_d^{a}$ denote set of the holomorphic foliations of degree $d\geq2$ admitting a compact curve as an integral curve.  Jouanolou \cite{J1979} showed that $\Bc_d^{a}$ is a proper subset of  $\Bc_d$. It follows that it is a countable union of proper analytic subsets in $\Bc_d$. For the higher dimension case, see \cite {LNSo1996}.

Let $\Hc_d$ denote the set of foliations of degree $d$ with only hyperbolic singularities. This is an open real Zariski dense set in $\Bc_d$. It is easy to show that $\Bc_d^{a}\cap\Hc_d$ is an analytic set  \cite{FS3} because one can bound the degree of the curve, when all singularities are hyperbolic.
Let $\Bc_d^{c}$ denote the set of foliations which admit directed positive closed currents. It is clearly a closed set in $\Bc_d$.
Part of Theorem 3.1 can be stated as $ \Bc_d^{c}\cap\Hc_d=\Bc_d^{a}\cap\Hc_d$.

 \begin{problem}
 Is the space $\Bc_d^{c}$ a countable union of proper analytic subsets of $\Bc_d$, when $d\geq2$ ?
 \end{problem}
 
Indeed, the only positive directed closed currents near a hyperbolic singularity or a saddle node singularity are supported by analytic sets.  It is tempting to think that $\Bc_d^c$ may be the set of points of positive Lelong number for some current in $\Bc_d$. The later is a countable union of analytic sets by a classical theorem of Siu, see \cite{Siu}.

\section{Anosov problem on the topology of leaves} \label{s:Anosov}

Consider foliations on $\P^2$. Recall that  $\Bc_d$ denote the space of holomorphic foliations of geometric degree $d$ on $\P^2$. Let  $\Ac_d$ denote the space of foliations given in an affine coordinate chart $\C^2$ in $\P^2$ by a vector field
$$Z:=P(z,w) {\partial \over \partial z} + Q(z,w) {\partial \over \partial w}$$
with $\max(\deg P, \deg Q)=d$. This notion of degree is not invariant under automorphisms of $\P^2.$
Generically in $\Bc_d$, the foliations have no compact leaf, and generically in $\Ac_d$, the line at infinity is invariant. 

\begin{problem}[Anosov]
For a generic foliation in $\Ac_d$ (resp. $\Bc_d$), are generic leaves biholomorphic to the unit disc ?
\end{problem}

Here, by generic leaves we mean all leaves except  countably many of them. The problem was discussed by Ilyashenko \cite{Il}. There are partial results for $\Bc_d$ by Sibony-Wold in \cite{SW}. For the case of $\Ac_d$, the question was discussed by Shcherbakov, Rosales-Gonz\'alez, Ortiz-Bobadilla in \cite{SRO}.

\section{Exceptional minimal set problem} \label{s:min}

Let $\Fc$ be a holomorphic foliation on $\P^2$. Denote by $\sing(\Fc)$ the singular set of $\Fc$ which is finite. 

\begin{problem}
Is it true that for every leaf $L$ of $\Fc$, we have $\overline L\cap \Sing(\Fc)\not=\varnothing$ ?
\end{problem}

The question is quite old and is due to Camacho-Lins Neto-Sad \cite{CLS}. The problem is also discussed in \cite{FS3}.  The answer is positive if $\Fc$ admits a directed positive closed $(1,1)$-current $T$, in particular, when $\Fc$ has an invariant algebraic curve. The support of this current $T$ should contain a singularity since otherwise the cohomology class of $T$ satisfies $\{T\}^2=0$ which is impossible in $\P^2$. Using the leaf $L$, we can construct a directed positive $\ddc$-closed current $S$ supported by $\overline L$. 
Again, for cohomological reason, one can show that $S$ and $T$ have to intersect and the intersection should be at the singularities as both currents are directed by the foliation. For this purpose, one can use the theory of densities for currents, see also \cite{DS1}.

\begin{problem}
Does there exist a Levi flat smooth or real analytic  hypersurface in $\P^2$ ? 
\end{problem}

Recall that a smooth or real analytic hypersurface $M$ of $\P^2$ can be defined as the zero set of a smooth or real analytic function $r$ in a neighbourhood of $M$ such that the differential $dr$ nowhere vanishes on $M$. We say that $M$ is {\it Levi flat} if in local coordinates $z=(z_1,z_2)$ we have 
$$\sum_{1\leq j,l\leq 2} {\partial^2 r(z) \over \partial z_j\partial \overline z_l} \xi_j\overline\xi_l =0 \ \text{ for every } \ \xi=(\xi_1,\xi_2)\in T_z^{1,0}(M).$$
Note that a tangent $(1,0)$-vector $\xi$ of $\P^2$ at $z\in M$ is in $T_z^{1,0}(M)$ if $\langle \partial r(z),\xi\rangle =0$. The above hypersurface $M$ is Levi flat if and only if it is foliated by Riemann surfaces.

It is shown in \cite[Example 7]{FS3} that no such hypersurface exists in $\P^2$ if we assume it is smooth and real algebraic. 
The proof is based on a Segre type argument as follows. 
Let $P(z,\overline z)=0$ be the equation of such a hypersuface in an affine chart of $\P^2$, with $P(z,\overline z)$ a real-valued polynomial in $z,\overline z$. Consider a germ of a leaf $L$ parametrized by a bi-holomorphic map $u:\D\to L$. Then, the function $P(u(t),\overline u(s))$ on $\D^2$ vanishes when $s=\overline t$, where $\overline u(s):=\overline{u(\overline s)}$. Since this function is a holomorphic, its zero set is a complex analytic set. Hence, it vanishes everywhere in $\D^2$. Fixing an $s$, we deduce that
$P(z, u(s))=0$ for $z\in L$. This is precisely the Segre argument which implies that the leaves are closed algebraic curves.  This is impossible because algebraic curves in $\P^2$ should intersect. 

\section{Function theory on foliations by Riemann surfaces} \label{s:RS}

The question is to study the function theory of leaves.
According to Brunella \cite{Br}, when a foliation $\Fc$ by Riemann surfaces in $\P^k$ has only hyperbolic singularities and no leaf supported by a compact curve, then there is no directed positive closed current. 

According to \cite{DS1}, when a foliation $\Fc$ in $\P^2$ has only hyperbolic singularities and admits leaves supported by compact curves, the only $\ddc$-closed currents directed by the foliation are the currents of integration on these compact curves.

A Riemann surface $Y$ is {\it parabolic} if all bounded subharmonic functions on it are constant. If a leaf of a foliation is parabolic, then the foliation admits a directed positive closed current of bi-dimension $(1,1)$, see \cite{BurnsSibony}. It follows that when there is no such current as in the case described by Brunella, the leaves are not parabolic and hence admit Green functions.

\begin{problem} (i) Assume that $\Fc$  is a foliation of geometric degree $d\geq 2$ on $\P^2$ with the line at infinity invariant and having only hyperbolic singularities. Do the leaves, except for the compact curves, have Green functions ?

(ii) More generally, assume that $\Fc$ is in $\Bc_d$ with $d\geq 2.$ What is the function theory of leaves (in terms of Green functions), when the foliation admits a non-zero, positive closed current directed by $\Fc$ ?
\end{problem} 

\section{Foliations with many directed harmonic currents} \label{s:many_currents}

Furstenberg \cite{Fu} has constructed uncountably many ergodic probability measures for the following map on the real torus $\T^2$ of dimension 2
$$T_\alpha(\xi_1,\xi_2) :=(e^{2i\pi\alpha}\xi_1, g_\alpha(\xi_1)\xi_2),$$
where $\alpha$ is a irrational number and $g_\alpha$ is a well chosen smooth function on the unit circle $\T$. 
In particular, every orbit of $T_\alpha$ is dense.

Based on his construction, Forn\ae ss-Sibony-Wold  constructed a compact minimal lamination $\Fc$ which admits uncountably many extremal directed closed currents which are mutually singular, see \cite{FSW}. Here, by minimality, we  mean that all leaves of $\Fc$ are dense. This lamination is defined on a compact Cauchy-Riemann manifold.

\begin{problem}
Construct holomorphic minimal foliation (resp. lamination) on a compact complex manifold with uncountably many extremal $\ddc$-closed, non-closed, currents
directed by the foliation (resp. lamination).
\end{problem}

\section{Transcendental vector fields} \label{s:trans}

The question is to explore the analogues of the results obtained for polynomial vector fields in $\C^k$ in the context of general holomorphic vector fields.

\begin{problem}
Let $Z$ be a holomorphic vector field in $\C^2$. 

(i) Is it possible to have a set of positive measure of leaves covered by $\C$ and another set of positive measure of leaves covered by the unit disc ?

(ii) If a non-pluripolar set of leaves are parabolic, are all leaves parabolic ?

(iii) Find conditions which imply that all leaves are covered by the unit disc.
\end{problem}

When the vector field is a polynomial vector field and all the singularities in $\P^2$ are non-degenerate, then all leaves are covered by the unit disc, see \cite {L}. For transcendental vector fields, the first "trivial" example is the one where the leaves are the level sets of a fixed entire transcendental function $f$. 
When $f$ is a submersion, this has been explored by Nishino and Yamaguchi in \cite {N, Y}.

\begin{problem}
Let $Z$ be a holomorphic vector field in $\C^k$ which is not the product of a polynomial vector field with a holomorphic function.
Let $\phi$ be the uniformization of a leaf $L$. The source could be $\C$ or the unit disc. Consider the Nevanlinna currents 
$$\tau^\phi_r:={1\over T(r)} \phi_*\Big(\log^+{r\over |\zeta|}\Big).$$
Are the cluster points of $\tau_r^\phi$, when $r$ tends to $\infty$ for the case of $\C$ and $r$ tends to $1$ for the case of the unit disc, supported necessarily by the hyperplane at infinity ?  
\end{problem}

First, one should observe as follows from \cite[Proposition 2.5]{DS1} that $T(r) \to\infty$ when
 $r$ tends to $1$.
In the case where a leaf $L=\phi(\D)$ is a closed, non-algebraic, analytic set in $\C^k$, the statement is true. Indeed, it follows from \cite[Proposition 2.6]{DS1}  that if there is a $\ddc$-closed $(1,1)$-currents on $\P^2,$  giving mass only to $\phi(\D)$, then it should be closed. Hence, the volume of $\phi(\D)$ is bounded in $\P^k$ and $\phi(\D)$ is algebraic.

\section{Dynamics on the space of maps} \label{s:maps}

In this section, we discuss a problem which is closely related to foliation theory and Nevanlinna's theory.
Let $\Hc$ denote the space of entire maps from $\C$ to $\C$ (resp. to $\P^1$), endowed with the topology of uniform convergence on compact sets. Let $\Uc$ be the set of $h\in\Hc$ such that for any other function $f\in\Hc$, there is a sequence $a_n\in \C$ such that $h_{a_n}(z):=h(z+a_n)$ converges to $f$.
G. D. Birkhoff \cite{B} has shown that $\Uc$ is non empty. 

Since $\Hc$ is a complete metric space with a countable basis, it follows that $\Uc$ is a dense $G_\delta$ set. The space
 $\Hc$ is laminated by images of $\C$ which are obtained by the translates of a given map. So it has  a $G_\delta$ dense set of dense orbits.
 The set of Brody maps (entire maps with bounded derivative with values in $\P^1$ and $\P^1$ is endowed with the Fubini-Study metric) is a dense set. If we fix a bound for the derivative, then we get a closed set of Brody maps. So Brody maps form an $F_\sigma$ set. Clearly, it does not intersect  $\Uc$, hence it has empty interior.

\begin{problem}
Find the minimal growth of the Nevanlinna characteristic function $T_h$ for $h\in\Uc$, where 
$$T_h(r):=\int_0^r {dt \over t} \int_{\D(t)} h^*(\omega)$$
and $\omega$ is the Fubini-Study form on $\P^1$. 
\end{problem}

Recently, Benjamin Weiss has proved the existence of non-trivial probability measures on $\Hc$ which are invariant by the above translation operator, see \cite{Weiss}. These measures are studied in particular by Buhovsky, Gl\"ucksam, Logunov and Sodin \cite{BGLS}. A measure is {\it non-trivial} if it gives zero measure to the set of constant maps. Using standard results from ergodic theory, one gets the existence of non-trivial
ergodic measures, i.e. measures such that the only invariant measurable subsets are of zero or full measure.

\begin{problem}
Given such a measure $\mu$ on $\Hc$ as above, what we can say about $\mu(\Uc)$ ?
\end{problem}

It is possible to consider the Birkhoff result to many situations : for the space $\Dc$ of holomorphic maps from the unit disc $\D$ to $\P^1$. Here, the translations on $\C$ are replaced by the holomorphic automorphisms of  $\D$. One can replace the unit disc by the unit ball or unit polydisc in
$\C^n$ and consider the maps with values in $\C^n$. The paper \cite{Be} contains a number of references.

\begin{problem}
Consider a complex manifold $M$ of dimension $n$ which is homogeneous under the action of a  Lie group $G$. Let $\Mc$ denote the space of holomorphic maps from $M$ to $\P^n$, endowed with the topology of uniform convergence on compact sets. The group $G$ acts on this space. 

(i) Under which conditions this action admits a dense orbit ?

(ii) Extend the construction by B. Weiss to this context.
\end{problem}

\begin{problem}
Explore the above properties for the following classes of maps in $\C^n$ : biholomorphisms of $\C^n$ and volume preserving (resp. symplectic) automorphisms of $\C^n$ (resp. $\C^n$ with $n$ even).
\end{problem}

\section{Ergodic theory for foliations with leaves of higher dimension} \label{s:HD}

Let $\Fc$ be a holomorphic foliation, possibly with singularities, in a compact K\"ahler manifold $M$, endowed with a fixed K\"ahler form $\omega$. When the leaves are Riemann surfaces, there is always a positive $\ddc$-closed current of bi-dimension $(1,1)$ directed by the foliation \cite{BS}. Indeed, the existence of such currents does not require any integrability condition, see \cite{Si}.

It is observed in \cite{FSW} that in general, the product of two foliations by Riemann surfaces does not admit a positive $\ddc$-closed current of bi-dimension $(2,2)$ which is directed by the product foliation. This is due to the fact that the product of two harmonic functions is not pluriharmonic in general. 

\begin{problem}
In order to study foliations with leaves of dimension larger than $1$, what is the right notion which substitutes directed $\ddc$-closed currents in dimension $1$ that should be considered ? 
Such objects should permit to develop the ergodic theory for these foliations.
\end{problem} 

\begin{problem}
Construct non-elementary  foliations $\Fc$ in $\P^k$ with leaves of dimension larger than $1$  and such that every leaf does not contain a compact curve. 
  \end{problem}
  
 We consider as {\it elementary} the foliations defined by closed forms or forms that become closed after multiplication by some integrating factors. For example, in the case of foliations by hypersurfaces, we disregard the situation where the foliation is given by a  logarithmic form with first integrals, i.e. the leaves are level sets of multivalued holomorphic functions of type  $f=f_1^{\lambda_1}\ldots f_p^{\lambda_p}$, with the $\lambda_j$'s real. In which case, the leaves could be dense and do not contain compact curves. 
  
The question is open even for foliations of codimension 1. Cerveau and Lins-Neto have described the components of the space of foliations of codimension 1 and of degree 2 in $\P^n,$ see \cite {CL}.

\section{Approximation of foliations with leaves of higher dimension} \label{s:app}

Runge approximation theorem is a central result in the theory of one and several complex variables. Consider now an analog in the theory of foliations. 

Let $\Fc$ be a holomorphic foliation in the unit polydisc $\D^k$ in $\C^k$. Assume that the leaves of $\Fc$ are of dimension $d$ with $d>1$. More precisely, the leaves of $\Fc$ are tangent to a family of $k-d$ holomorphic $(1,0)$-forms  $(\alpha_j)_{1\leq j\leq k-d}$ on $\C^k$. These forms $\alpha_j$ are linearly independent at a generic point of $\D^k$  and satisfy the Frobenius condition 
$$\partial \alpha_j =\sum_{l=1}^{k-d} \alpha_l \wedge \beta_{j,l} \quad \text{with} \quad j=1,\ldots,k-d$$
for some holomorphic $(1,0)$-forms $\beta_{j,l}$ in $\D^k$. The last condition is equivalent to 
$$\partial \alpha_j\wedge \alpha_1\wedge \ldots\wedge \alpha_{k-d} =0 \quad \text{with} \quad j=1,\ldots,k-d.$$

\begin{problem}
Is it always possible to approximate the foliation $\Fc$ by foliations in $\C^k$ or\break $\P^k$ ? More precisely, is it always possible to approximate the forms $\alpha_j$, uniformly on compact sets of $\D^k$, by entire forms (resp. polynomial forms) $\alpha_{j,n}$ which still satisfy the above Frobenius  condition ?
\end{problem}

The question appears in a problem of Bogomolov \cite {Bo} for foliations in a Stein manifold. He conjectures that any compact complex manifold admits an embedding as a smooth subvariety transverse to an algebraic  foliation on a projective algebraic manifold.

\end{document}